\documentclass[11pt,psamsfonts]{amsart}

\usepackage{amsmath}
\usepackage{amsthm}
\usepackage{amssymb}
\usepackage{amscd}
\usepackage{amsfonts}
\usepackage{amsbsy}
\usepackage{epsfig,afterpage}
\usepackage{psfrag}

\usepackage{graphicx}

\newtheorem {theorem} {Theorem}

\newtheorem {remark} [theorem]{Remark}

\begin{document}

\title[Periodic orbits in the hyperchaotic Chen system]
{Periodic orbits in the hyperchaotic Chen system}

\author[S. Maza]
{ Susanna Maza$^1$}

\address{\noindent $^1$ Departament de Matem\`atica, Universitat de
Lleida, Avda. Jaume II, 69, 25001 Lleida, Catalonia, Spain}

\thanks{The author is partially supported by a MICINN grant number MTM2011-22877 and by a CIRIT grant number 2014 SGR 1204.}



\begin{abstract}
In this work, we show that a zero--Hopf bifurcation take place in the Hyperchaotic Chen system as parameters vary. Using averaging theory, we prove the existence of two periodic orbits bifurcating from the zero--Hopf equilibria located at the origin of the Hyperchaotic Chen system.
\end{abstract}

\maketitle

\section{Introduction and statement of the results}\label{s1}

Hyperchaos has been widely investigated since in 1979 Otto R\"ossler proposed one of the first hyperchaotic attractors, see \cite{Ross}. A hyperchaotic system is a chaotic one in which two or more Lyapunov exponents are positive indicating that the chaotic dynamics are expanded in more than one direction. This feature gives rise to complex attractors harder to control than chaotic ones. Thus, hyperchaos has been found to be useful in many fields such as in encryption, secure communications \cite{chen2} and also is matter of interest in nonlinear circuits \cite{li}, liquid mixing, lasers \cite{hak} and many more.

It is worth to say that the minimal dimension of the phase space that embeds a hyperchaotic attractor must be four, so the typical examples of hyperchaotic systems have been introduced usually as extensions of known autonomous three--dimensional chaotic systems. Besides the hyperchaotic attractor of R\"ossler, the  hyperchaotic Lorenz-Haken system \cite{ning}, the Chua's circuit \cite{chua} and the hyperchaotic L\"u system \cite{lu} are well--known examples of hyperchaotic models as extensions of three--dimensional chaotic ones. The book \cite{Sprott} is a catalog of systems exhibiting chaos.
\smallskip

In this work we are interested in the hyperchaotic Chen system,
\begin{equation}\label{chen-1}
\begin{array}{ll}
\dot{x} =  a (y - x) + w, \\
\dot{y} = d x + c y - x z,  \\
\dot{z} = x y - b z,\\
\dot{w} =y z + r w.
\end{array}
\end{equation}
System (\ref{chen-1}) was introduced in \cite{chen} also as an extension of the three--di\-men\-sio\-nal chaotic Chen system. Its study has generated considerable research, techniques of chaos synchronization, studies of controlling chaos, secure communications, power system protection and so on, see \cite{he}, \cite{effa}, \cite{park}, \cite{smaoui} and the references therein.

The dynamics of Hyperchaotic systems has been studied mainly numerically and by computer simulations. In this work we perform an analytic analysis on the existence of periodic orbits of the differential system (\ref{chen-1}) by applying {\it averaging theory of first order}. More precisely, we will prove that a {\it zero--Hopf} bifurcation occurs in system (\ref{chen-1}) bifurcating two limit cycles from the {\it zero--Hopf} equilibria as parameters vary. As far as I know, this is the first time that an analytic analysis of a zero--Hopf bifurcation in hyperchaotic Chen system is performed.

We recall that an equilibrium of a differential system $\bf{\dot{x}}=\bf{f(\bf{x})}$  with ${\bf f}:A \to \mathbb{R}^n$ and $A $ an open subset of $\mathbb{R}^n$  is a zero--Hopf equilibrium if it has two pure imaginary conjugated eigenvalues and $n-2$ zero eigenvalues. The bifurcation of periodic orbits from zero--Hopf equilibria of three dimensional differentials systems has been studied via averaging theory in \cite{lli2}. See also \cite{glm2} for a study in systems with two slow and one fast variables. In \cite{glm} is performed an analysis of the periodic orbits bifurcating in the chaotic prototype-4 system of  R\"ossler. Recently,  was  published \cite{lli1}, a paper on a zero--Hopf bifurcation in the hyperchaotic lorenz system.
\smallskip

The main result of this work is the following.

\begin{theorem}\label{Teo-chen}
Consider the hyperchaotic  Chen system {\rm (\ref{chen-1})} with $c=a$ for small values of $b$ and $r$, and satisfying $b(a + d) \ r>0$ and $a(a+d)<0$. Then a zero--Hopf bifurcation occurs emerging two limit cycles from the origin.
\end{theorem}

In section 2 we give a summary of the averaging theory of first order for finding limit cycles and we proof Theorem \ref{Teo-chen} in section 3.

\section{Averaging theory: perturbing an isochronous system}

We consider the problem of bifurcation of $T$--periodic solutions for a differential system of the form
\begin{equation}\label{pert-A1}
\dot{{\bf x}}= {\bf F_0}(t,{\bf x}) + \varepsilon {\bf F_1}(t,{\bf x}, \varepsilon) \ ,
\end{equation}
where $\varepsilon$ is a small positive parameter, ${\bf x} \in A$, where $A$ is an open subset of $\mathbb{R}^n$, $t\geq 0$. Moreover we assume that both ${\bf F_0}(t,{\bf x})$ and ${\bf F_1}(t,{\bf x}, \varepsilon)$ are
$\mathcal{C}^2$ functions and $T-$periodic in $t$.

The classical theory of averaging reduce the problem of finding $T$--periodic solutions of (\ref{pert-A1}) for $\varepsilon>0$  small to the problem of finding simple zeros of the so--called {\it bifurcation functions}. Many methods encountered in the literature are based on this idea, see \cite{buica},\cite{ma}, \cite{chicone3}.
\newline
We study here the particular case in which all the solutions of the unperturbed system
\begin{equation}\label{unper}
 \dot{{\bf x}}= {\bf F_0}(t,{\bf x})
\end{equation}
are $T$--periodic. Then we assume that the unperturbed system is i\-so\-chro\-nous and the problem in what we are interested in is the problem of the persistence of the periodic orbits under some perturbation.

We denote the linearization of (\ref{unper}) along a periodic solution ${\bf x}(t,{\bf u})$  of (\ref{unper}) such that ${\bf x}(0,{\bf u})={\bf u}$ by

$$
\dot{{\bf y}}=D_{\bf x}{\bf F_0} (t,{\bf x}(t,{\bf u}))
$$
where $D_{\bf x}{\bf F_0} $ is the Jacobian matrix of $F_0$ with respect to ${\bf x}$ and let $\Phi_{\bf u}(t)$ be some fundamental matrix of (\ref{unper}).
Assume that there exists an open set $V$ with closure $\bar{V} \subset A$ such that for each ${\bf u}\in \bar{V}$, ${\bf x}(t,{\bf u})$ is $T$--periodic. The following result gives an answer to the question of bifurcating periodic solutions
from the $T$--periodic solutions ${\bf x}(t,{\bf u})$.

\begin{theorem}\label{ave}
We assume that there exists an open set $V$ with $\bar{V} \subset A$ such that for each ${\bf u} \in \bar{V}$, ${\bf x}(t,{\bf u})$ is $T$--periodic. Consider the function ${\bf f}: \bar{V}\mapsto \mathbb{R}^n $ given by
\begin{equation}\label{bifur}
 {\bf f}({\bf u})= \frac{1}{T} \int_0^{T}\Phi_{\bf u}^{-1} {\bf F_1}(t,{\bf x}(t,{\bf u})) \, dt \ . \\
\end{equation}
If there exists ${\bf p} \in V$ with ${\bf f}({\bf p})=0$ and $\det(D_{\bf u}{\bf f}({\bf  p})))\neq0$ then there exists a T--periodic solution  $\gamma(t,\varepsilon)$ of system {(\rm \ref{pert-A1})} such that  $\gamma(0,\varepsilon)\to {\bf  p}$ as $\varepsilon \to 0$. Moreover, if all the eigenvalues of $\det(D_{\bf u}{\bf f}({\bf  p})))$  have negative real part, the corresponding periodic orbit $\gamma(t,\varepsilon)$ is asymptotically stable for $\varepsilon$ sufficiently small.
\end{theorem}

For a proof of Theorem \ref{ave} see \cite{ma}, \cite{Rose} and \cite{buica2}.

\section{Proof of Theorem \ref{Teo-chen}}

The origin of system (\ref{chen-1}) is always an equilibrium point of it under under any choice of parameters. The characteristic polynomial $p(\lambda)$ of the linearization of system (\ref{chen-1}) at the equilibrium point located at the origin is given by
$$
p(\lambda) = (r - \lambda)(b + \lambda)(a (c + d - \lambda) + (c - \lambda)\lambda)
$$
The eigenvalues associated at the origin are
$$ \lambda_1=r, \ \lambda_2=-b, \ \lambda_{3,4}= \frac{1}{2} \left(-a + c \pm\sqrt{a^2 + 2ac + c^2 + 4ad} \right).
$$
This suggest to consider small values of the parameters $r$ and $b$ introducing the small parameter $\varepsilon$ in the following way $(r,b)\mapsto (\varepsilon r,\varepsilon b)$. Thus, taking $c=a$ the eigenvalues of system (\ref{chen-1}) at the origin becomes $\lambda_{1,2}=0$ and $\lambda_{3,4}=\pm \sqrt{a(a+c)}$ when $\varepsilon\mapsto 0$. It follows that the origin is a zero--Hopf equilibrium when $\varepsilon\mapsto 0$ if $a(a+c)<0$. So, for this choice of parameters we have system (\ref{chen-1}) written as
\begin{equation}\label{chen-2}
\begin{array}{ll}
\dot{x} =  a (y - x) + w, \\
\dot{y} = d x + a y - x z,  \\
\dot{z} = x y - b\varepsilon  z,\\
\dot{w} =y z + r\varepsilon w.
\end{array}
\end{equation}
Reescaling the variables $(x,y,z,w)\mapsto (\varepsilon x,\varepsilon y,\varepsilon z,\varepsilon w)$, system (\ref{chen-2}) becomes
\begin{equation} \label{chen-3}
\left( \begin{array}{c} \dot{x} \\ \dot{y}\\ \dot{z} \\ \dot{w}  \\
\end{array} \right) =\left(
\begin{array}{c}  a (y - x) + w \\ d x + a y \\0 \\0 \\
\end{array} \right) \, + \varepsilon  \left(\begin{array}{c} 0 \\- x z \\ x y - b z \\ y z + r w \end{array} \right).
\end{equation}
Now we have Chen system (\ref{chen-1}) written as differential system of the form (\ref{pert-A1}) and we can consider the problem of bifurcating periodic solutions of it by using averaging theory.
We have to solve first the unperturbed system of (\ref{chen-3}). The solution ${\bf x}(t,{\bf u})= (x(t),y(t),z(t),w(t))$ of
\begin{equation}\label{chen-4}
\begin{array}{ll}
\dot{x} =  a (y - x) + w, \\
\dot{y} = d x + a y,   \\
\dot{z} = 0, \\
\dot{w} =0, \\
\end{array}
\end{equation}
satisfying the initial condition $ {\bf u}=(x(0),y(0),z(0),w(0))=(x_0,y_0,z_0,w_0) \in \mathbb{R}^4$ is
\begin{eqnarray}
x(t) &=& \frac{1}{a + d} (w_0 +((a + d)x_0-w_0)\cosh{(\sqrt{a(a + d)} \ t)}+ \nonumber \\
 & & \frac{\sqrt{a + d}}{\sqrt{a}}(w_0 + a(y_0-x_0))\sinh{(\sqrt{a(a + d)}\ t})), \nonumber \\
y(t) &=& \frac{1}{a(a + d)}((d w_0 + a(a + d)y_0)\cosh{(\sqrt{a (a + d)} \ t)}-d w_0+  \label{chen-5} \\
 & & \sqrt{a}\sqrt{a + d}(dx_0 + ay_0)\sinh{(\sqrt{a(a + d) }\ t)}), \nonumber \\
w(t) &=& w_0, \nonumber \\
z(t) &=& z_0. \nonumber
\end{eqnarray}
Notice that if $a (a + d)<0$ then $\cosh{(\sqrt{a(a + d)} \ t)}=\cos{(\sqrt{-a(a + d)} \ t)}$ and $\sinh{(\sqrt{a(a + d)} \ t)}= i \sin{(\sqrt{-a(a + d)}\ t)}$
being $i^2=-1$. Hence, when  $a(a + d)<0$ the solution (\ref{chen-5}) of the unperturbed system (\ref{chen-4}) is
\begin{eqnarray}
 x(t) &=& \frac{1}{a + d} (w_0 +((a + d)x_0-w_0)\cos{(\Omega t)}- \nonumber \\
 & & \frac{\Omega}{a}(w_0 + a(y_0-x_0))\sin{(\Omega \ t})), \nonumber \\
y(t) &=& \frac{1}{a(a + d)}((d w_0 + a(a + d)y_0)\cos{(\Omega t)}- \label{chen-6} \\
 & & d w_0- \Omega(dx_0 + ay_0)\sin{(\Omega t)}), \nonumber \\
w(t) &=& w_0, \nonumber \\
z(t) &=& z_0. \nonumber
\end{eqnarray}
where $\Omega=\sqrt{-a(a + d)}$. We have that any solution (\ref{chen-6}) of the unperturbed system (\ref{chen-4}) is periodic of period $T=\frac{2 \pi}{\Omega}$. So system (\ref{chen-3}) is, in fact, a perturbation of an isochronous system when $a(a+c)<0$ and we can apply Theorem \ref{ave}.
The first variational system of (\ref{chen-4}) along the solution (\ref{chen-6}) coincides with the unperturbed system (\ref{chen-4}), so the inverse of fundamental matrix solution $\Phi_{\bf u}(t)$ of (\ref{chen-4}) is {\small
$$
 \left(\begin{array}{cccc} \cos{(\Omega t)}+ \frac{a}{\Omega}\sin{(\Omega t)}& -\frac{a}{\Omega}\sin{(\Omega t)} & 0 & \frac{1}{(a+d)} (1-\cos{(\Omega t)}+ \frac{\Omega}{a} \sin{(\Omega t)})  \\ -\frac{d}{\Omega}\sin{(\Omega t)}  & \cos{(\Omega t)}- \frac{a}{\Omega}\sin{(\Omega t)} & 0 & \frac{d}{a (a+d)}(\cos{(\Omega t)})-1) \\ 0 & 0 & 1 & 0 \\ 0 & 0 & 0 & 1 \end{array} \right).
$$}
The bifurcation function (\ref{bifur}) is given by
 	
$$
 f({\bf u})= \frac{\Omega}{2 \pi} \int_0^{\frac{2 \pi}{\Omega}}\Phi_{\bf u}^{-1} F_1(t,{\bf x}_{\bf u}) \, dt = (f_1,f_2,f_3,f_4) \ , \\
$$
being
\begin{eqnarray*}
f_1({\bf u}) &=& \frac{r w_0}{a+d} + \frac{dz_0((a + d)x_0-3w_0)}{2a(a + d)^2}- \frac{(a(y_0 - x_0)+w_0)z_0}{2(a + d)}\ , \\
f_2({\bf u}) &=& \frac{d(3dw_0 + a(a + d)y_0)z_0} {2 a^2 (a + d)^2}- \frac{d r w_0}{ a (a + d)}-\frac{(d x_0+ a y_0)z_0 } {2 (a + d)} , \\
f_3({\bf u}) &=&  \frac{ d^2 (x_0^2 + 2 b z_0)}{2 (a + d)^2}+ \frac{a(2abz_0-y_0(2w_0 + a(y_0-2x_0))}{2 (a + d)^2}- \\
& & \frac{(3w_0^2+2aw_0y_0) d}{2 a(a + d)^2} +\frac{a d(x_0^2 + 2x_0y_0 - y_0^2 + 4bz_0)) }{2 (a + d)^2}  \ , \\
f_3({\bf u}) &=& w_0 \left( r - \frac{d z_0}{a (a + d)}\right) \ .
\end{eqnarray*}

The isolated zeros of the map $ {\bf u} \mapsto f({\bf u}) = (f_1({\bf u}),f_2({\bf u}),f_3({\bf u}),f_4({\bf u}))$  are
\begin{eqnarray*}
 {\bf p}_{1,2}&=& \left(\pm\frac{a \sqrt{b(a + d) \ r}}{d} ,\mp \sqrt{ b(a + d) \ r} ,\frac{a (a + d) \  r}{d},\pm \frac{a \sqrt{b (a + d)\ r}\ (a + d)}{d} \right) \\
 \end{eqnarray*}
Notice that ${\bf p}_{1,2} \in \mathbb{R}$ because of conditions of Theorem \ref{Teo-chen}.
The determinant of the Jacobian matrix of ${\bf f}$ at the points ${\bf p}_{1,2}$ is
$$
\det(D {\bf f}({\bf p_{1,2}})) = \frac{b (a^4 + a^3 d - d^2) r^3}{2 d^2}.
$$
From condition $a(a+d)<0$ we have that $d\neq0$ and  $a^4 + a^3 d - d^2= a^2(a(a+d)-d^2)\neq0$. Moreover, from  $b(a+d)r<0$ we have $b r\neq0$. Then we get $\det(D {\bf f}({\bf p_{1,2}}))\neq0$ and ${\bf p}_{1,2} $ are simples zeroes of ${\bf f}$. Hence, the averaging theory stated in Theorem \ref{ave} predicts the existence of two $T$--periodics orbits $\gamma_{1,2}(t,\varepsilon)$ of system (\ref{chen-2}) with period $\frac{2 \pi}{\Omega}$  such that $\gamma_{1,2}(0,\varepsilon)\to (p_{1,2})$ as $\varepsilon \to 0$.
\newline

Since we have performed the reescaling $(x,y,z,w)\mapsto (\varepsilon x,\varepsilon y,\varepsilon z,\varepsilon w)$ for bringing system (\ref{chen-1}) to system (\ref{chen-2}), the solutions $\gamma_{1,2}(t,\varepsilon)$
of system (\ref{chen-2}) provides the periodic orbits $\varepsilon \gamma_{1,2}(t,\varepsilon)$ of system (\ref{chen-1}) tending to the zero-Hopf equilibrium as $\varepsilon \to 0$. This finishes the proof.

\begin{remark}
{\rm Regarding the stability of the bifurcated periodic orbits, we get that the eigenvalues of $D {\bf f}({\bf p_{1,2}})$ are
$$
\lambda_{1,2} =\frac{1}{2} \left(b \pm \sqrt{b (b + 8 r)} \right ), \ \   \lambda_{3,4}=\frac{r}{2} \left(1 \pm i \frac{a \Omega}{d}\right ) \ .
$$
Therefore, under our parameter restrictions we cannot have all the eigenvalues with negative real part. In consequence we cannot use Theorem \ref{ave} for the analysis of the periodic orbits $\gamma_{1,2}(t,\varepsilon)$.}
\end{remark}

\end{document}